\documentclass[11pt, a4paper, twocolumn]{article} 
\usepackage{graphicx}
\usepackage{subcaption}

\usepackage[english]{babel} 

\usepackage{microtype} 

\usepackage{amsmath,amsfonts,amsthm} 

\usepackage[svgnames]{xcolor} 

\usepackage{booktabs} 

\usepackage{lastpage} 

\usepackage{graphicx} 

\usepackage{enumitem} 
\setlist{noitemsep} 

\usepackage{sectsty} 
\allsectionsfont{\usefont{OT1}{phv}{b}{n}} 

\usepackage{geometry} 

\usepackage[T1]{fontenc} 
\usepackage[utf8]{inputenc} 

\usepackage{XCharter} 

\usepackage{fancyhdr} 
\fancypagestyle{firstpage}{ 
	\fancyhf{}
}

\newcommand{\institution}[1]{{\footnotesize\usefont{OT1}{phv}{m}{sl}\color{Black}#1}} 

\usepackage{titling} 

\newcommand{\HorRule}{\color{DarkGoldenrod}\rule{\linewidth}{1pt}} 

\pretitle{
	\vspace{-30pt} 
	\HorRule\vspace{10pt} 
	\fontsize{32}{36}\usefont{OT1}{phv}{b}{n}\selectfont 
	\color{DarkRed} 
}

\posttitle{\par\vskip 15pt} 

\preauthor{} 

\postauthor{ 
	\vspace{10pt} 
	\par\HorRule 
	\vspace{20pt} 
}

\usepackage{lettrine} 
\usepackage{fix-cm}	

\usepackage{xstring} 

\newcommand{\AuthorNames}{Eleuterio F.  Toro}
\newcommand{\ShortTitle}{Personal Reminiscences of S.  K.  Godunov}

\title{\sc Reminiscences of 
S. K.  Godunov\\The Russian  Mathematician} 

\author{\sc Eleuterio F.  Toro
	\newline\newline 
	 \institution{\sc Professor Emeritus,  Mathematics, University of Trento, Italy\\
	 Life Fellow, Clare Hall, University of Cambridge, UK\\
	 \tt Personal webpage: https://eleuteriotoro.com \\
	  Email: eleuterio.toro@unitn.it}}
	 
\date{} 

\usepackage{fancyhdr}
\begin{document}

\maketitle

\thispagestyle{firstpage}




\begin{center}
\begin{minipage}{0.95\columnwidth}
\noindent\textbf{\small{ABSTRACT}} 

{{These personal reminiscences of the great Russian mathematician Sergey K. Godunov (1929–2023) arose from a request by his daughter, Ekaterina, to contribute a piece to a book she is writing about her father’s life.  I was honoured to accept this invitation and to give written form to the rewarding experience of conducting research on themes pioneered by Professor Godunov, interacting with him personally on several memorable occasions, and helping to establish research collaboration with his Novosibirsk group.

Our association began at a conference in Lake Tahoe (USA) in 1995 and was followed by a number of subsequent meetings, notably in Novosibirsk, Manchester, Oxford, and Cambridge. Briefer encounters also took place in the Porquerolles Island (France),  in Lyon (France), and in St. Petersburg (Russia).

These notes bear witness to the global impact of Godunov’s mathematical creativity across multiple branches of science, as well as to its lasting influence on the careers of generations of mathematicians in both academia and industry.
}}
\end{minipage}
\end{center}



\begin{center}
{\bf Encountering Godunov's Method} 
\end{center}

As a researcher, over a considerable period of time I have followed the {\it Godunov School}, that is, the physically motivated manner of approaching the field of hyperbolic partial differential equations and the design of computational algorithms to solve them.

As a post-doctoral research fellow in the early 1980s at the  Cranfield Institute of Technology, UK, I vividly remember being thrown into the deep waters of compressible fluid dynamics, combustion, propulsion, shock waves, and hyperbolic partial differential equations. My previous PhD days working with finite element methods for elliptic free-boundary problems appeared far away—and even rosy.

Professor J. F. Clarke (1927–2013), at Cranfield,  had recently come across a paper on the Random Choice Method (RCM), or Glimm's method  \cite{Glimm:1965a}, applied to the Euler equations of gas dynamics and was impressed by its unrivalled resolution of shock waves and contact surfaces. He then encouraged me to try the wonders of RCM, but inevitably,  I also ended up coding Godunov’s method \cite{Godunov:1959a}.

\begin{figure}[h]
      \centerline{
             \includegraphics[scale=0.33, angle=0]{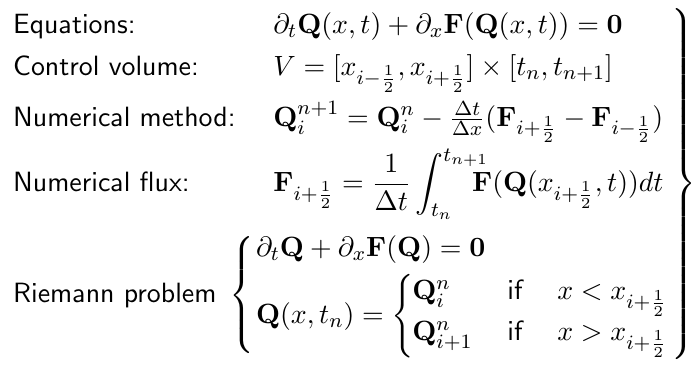}
      }
      \caption{{\bf Godunov's method. } {\it Equations} are integrated in {\it Control volume} $V$ leading to {\it Numerical method} to advance the solution from time $n$ to $n+1$. Crucially,  {\it Numerical flux} is determined by solving the {\it Riemann problem} with piece-wise constant data at time $t=t_{n}$; its solution is inserted in the time integral to determine the {\it Numerical flux},  and thus the {\it Numerical method}}
      \label{fig:GodunovMethod}
\end{figure}

Fig. \ref{fig:GodunovMethod} displays all key components of Godunov’s method.   {\bf Equations} symbolize any system of hyperbolic conservation laws,  potentially applicable to many areas of science and technology.  {\bf Control volume} $V$ is a representative element of a discretized computational domain in space-time. {\bf Numerical method} arises from exact integration of {\bf Equations} in control volume $V$, leading to integral averages for the unknown vector {\bf Q} and physical flux {\bf F}.  Formula {\bf Numerical method} serves to advance the solution in time  -- in the form of integral averages--from a known time level $n$ to the future unknown time level $n+1$.  However,  the formula {\bf Numerical method} requires the prescription of the {\bf Numerical flux} embodied in the time integral,  which in turn requires the unknown vector 
${\bf Q}$ at the interface as function of time $t$.  And here comes what is central to Godunov's method.  Godunov proposed to solve the Riemann problem, whose solution would give that needed function under the integral.   The Riemann problem is defined by two items,  the governing {\bf Equations} and local initial conditions at the known time level $n$ consisting of two constant vectors (integral averages),  with a discontinuity right at the interface.  In his paper of 1959 \cite{Godunov:1959a} Godunov proposed to solve the Riemann problem exactly and in 1962  \cite{Godunov:1962a} he presented a linearised approximation,  as an alternative.

At Cranfield, my first challenging task was to compute the exact solution of the Riemann problem for the Euler equations of gas dynamics. Somewhere, I found Godunov’s original solution method, the core of which was an iterative scheme for three non-linear algebraic equations. I remember the method as being exceedingly slow, requiring many iterations to converge.   Having a program to compute the exact solution of the Riemann problem, I implemented the solution in Glimm's method  \cite{Glimm:1965a}.  And yes, the resolution of shock waves and contact surfaces was impressive, though the smooth parts of the flow—somehow the easy part—were riddled with random noise.

With an exact Riemann solver available, I then ventured into territory that was, strictly speaking, outside my remit. I used the solution to compute the Godunov numerical flux; see Fig. \ref{fig:GodunovMethod}. My somewhat erratic efforts at the time ran counter to the trends of the 1980s, when approximate Riemann solvers were in fashion. Philip Roe had recently joined our department and led that trend with his newly published approximate, linearized Riemann solver \cite{Roe:1981a}.  Jack Pike, who had earlier worked with Roe on a Roe-Pike approximate Riemann solver,  also joined our group as a glorified post-doctoral fellow; the research environment was thrilling.

In evaluating the performance of Riemann solvers,  I found that the exact solver was rather expensive.  This was due to two factors. First, the exact solution itself is computationally more involved,  and second, Godunov’s particular three-equation iterative scheme was more complex than necessary.   I then developed an alternative iterative solver based on just one non-linear algebraic equation, which was published several years later \cite{Toro:1989d}.  In fact, it turns out that for the ideal Euler equations this exact Riemann solver is only marginally more expensive than any approximate Riemann solver; surprisingly,  it is cheaper than the Osher-Solomon approximate Riemann solver \cite{Osher:1982a},\cite{Toro:2009a}.  The situation is different for materials with complex equations of state.

\begin{figure}[h]
      \centerline{
             \includegraphics[scale=0.45, angle=0]{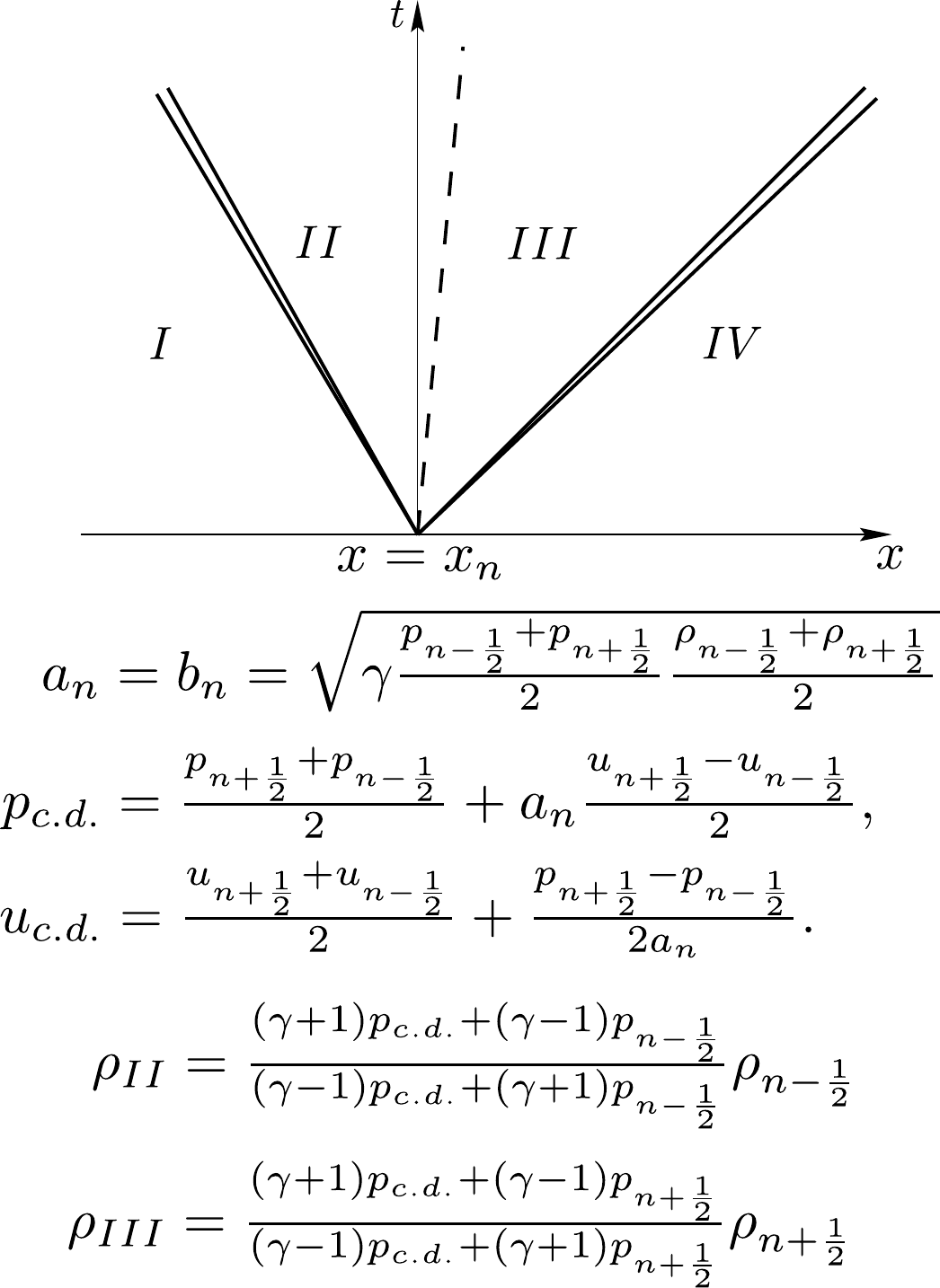}
      }
      \caption{Godunov's approximate Riemann solver for the Euler equations (1962) (reproduced from  the original publication  
      \cite{Godunov:1962a})}
      \label{fig:GodunovLinearised}
\end{figure}

In my beginner’s naivety,  I was tricked by the belief that Godunov’s method \cite{Godunov:1959a} was, by definition, the conservative finite-volume formula in which the numerical flux is determined by solving the Riemann problem exactly, perhaps using the old three-equation algorithm. For quite some time, I was unaware that Godunov himself had been the first to propose a linearized approximate scheme for solving the Riemann problem \cite{Godunov:1962a}.  Fig.  \ref{fig:GodunovLinearised} illustrates the wave pattern associated with the Riemann problem for the Euler equations, along with the explicit formulae for the solution in zones II and III, usually called the {\it Star Region} 
\cite{Toro:2009a}.  As a poor consolation, I realized that I was not alone in this misunderstanding. Even today, this misconception remains widespread, especially among young researchers. 

It is therefore appropriate to define a Godunov method as one that encompasses all the ingredients shown in Fig. \ref{fig:GodunovMethod}, in which the {\it Riemann problem} used to determine the time integral, and hence the {\it Numerical flux}, is solved either approximately or exactly.  In the last 40 years there have been significant advances in the developments of exact and approximate Riemann solvers in various fields of applications; see for example \cite{Toro:2009a, Toro:2024a,Toro:2025a} and references therein. 

Furthermore,  Fig.  \ref{fig:GodunovMethod} has provided a framework not only for the development of new variations of Godunov's original method but also for the interpretation of existing methods and the invention of novel ones.  For example,  and this may appear surprising to some, the original  Lax-Friedrichs method (or the Lax method) \cite{Lax:1954a} may be derived as the limiting case of an HLL-type Godunov method by assuming a fixed and symmetric wave pattern in Fig.  \ref{fig:GodunovLinearised},  with  the fastest wave speeds chosen as $\pm \Delta x/\Delta t$.   Moreover,  the Lax-Friedrichs updated state  ${\bf Q}_{i}^{n+1}$ in Fig.  \ref{fig:GodunovMethod} may be obtained from a spatial integral average at time $t=\frac{1}{2}\Delta t$ of the solution of the Riemann problem with data ${\bf Q}_{i-1}^{n}$ and ${\bf Q}_{i+1}^{n}$; for details see \cite{Toro:2009a},  chapters 5 and 7.

Further, the second-order,  non-monotone Lax-Wendroff method may be interpreted as a Godunov method by a suitable integration of the wave pattern arising from the solution of the Riemann problem to obtain the numerical flux \cite{Toro:1989a}.

The very-high order ADER methods developed in recent years, e.g. \cite{Toro:2001c, Toro:2002a, Titarev:2002a, Dumbser:2005a, Dumbser:2008a, Toro:2024b}, may also be interpreted as Godunov-type methods, in which the constant vectors in the Riemann problem of Fig.  \ref{fig:GodunovMethod} are replaced by functions of distance $x$, such as polynomials of degree $m$.   Moreover, suitable modifications of the framework of Fig.  \ref{fig:GodunovMethod} give rise to still larger families of numerical methods for solving partial differential equations.

\begin{center}
{\bf Lake Tahoe, USA} 
\end{center}

It was at a conference in Lake Tahoe, USA, in 1995 that I first caught sight of Sergey Godunov.  He was always surrounded by many people during the breaks; greeting him—and better still, talking to him—seemed challenging. In fact, I missed my chance on that occasion because, due to an unexpected health problem, I had to cut short my participation in the conference and return to Manchester.

My post-doctoral research fellow,  Dr.  Stephen Billett (1969-2022), stayed on and kindly agreed to present my paper. Back in Manchester, he told me that when he realized that Godunov himself was sitting in front of him, ready to listen to his talk, he went into a panic.  He was also aware that Russian scientists tended to stand up when asking a question,  which may prove rather intimidating to some.  As a matter of fact, Godunov did ask some questions,  standing up.  Stephen survived the experience and felt quite proud of it.

A couple of weeks later, I received a long, typed letter from Sergey Godunov. He was very complimentary about the work presented at the Lake Tahoe Conference.  Stephen had done a good job, obviously. The talk was concerned with anomalies of conservative methods in the presence of contact surfaces and shear waves \cite{Toro:2002b}. Then, a couple of years later, in 1997, I again had the chance to see Godunov from a distance, this time in a lecture theatre at the University of Michigan (USA). On that occasion, Godunov was awarded the distinction of {\it Doctor Honoris Causa}. But again, I missed the chance to meet the man.

\begin{center}
{\bf Novosibirsk,  Siberia, Russia}
\end{center}

It was in the beautiful Siberian summer of 1998, with temperatures reaching $28^\circ C$, that I first met Sergey Godunov in Novosibirsk. I participated in an all-Russian conference held from 22nd to 27th June, the Third Siberian Congress on Industrial and Applied Mathematics, dedicated to the memory of S. L. Sobolev (1908–1989). My contribution to the conference was titled: {\it Anomalies of total-energy conservative methods-- analysis, numerical evidence and possible cures}.

Due to language difficulties, I missed most of the talks. My old Russian course for mathematics students at the Warwick University Mathematics Institute, UK, was of no help. Sir Christopher Zeeman (1925–2016), as Head of the Institute, is reputed to have said that half of the world’s mathematics was produced in Russia, and therefore all mathematics students must learn Russian.

In compensation however, I had time to get to know the famous University of Novosibirsk. A mathematics PhD student kindly showed me around. He took me to the renowned Sobolev Institute of Mathematics and the also famous Institute of Hydrodynamics. The Sobolev Institute of Mathematics in Novosibirsk was founded in 1957 under the leadership of mathematicians Sobolev, Lavrentiev, and Khristianovich. The list of famous Novosibirsk mathematicians is long and includes illustrious names.  It is impressive to see the effigies of some of them crowning the building façade:
A. A.  Lyapunov,   L. V.  Kantorovich,  A.I.  Shirshov,  M.I.  Kargapolov,  A.I.  Maltsev,  S.L.  Sobolev and A.D.  Aleksandrov.
 
One evening, upon returning to my hotel, the concierge passed on a message reminding me that the following morning, at 11:00 o’clock, Professor Godunov would come to collect me at the hotel. That morning marked the first time we formally introduced ourselves and exchanged greetings. We also exchanged books on the occasion. Professor Godunov was known to be very proficient in French, less so in English, and my Russian was not up to scratch. Nonetheless, we managed to communicate quite well in English. He informed me that we would be walking to his dacha, some 6 to 7 km away in the countryside.
\begin{figure}[h]
      \centerline{
            \includegraphics[scale=0.53, angle=0]{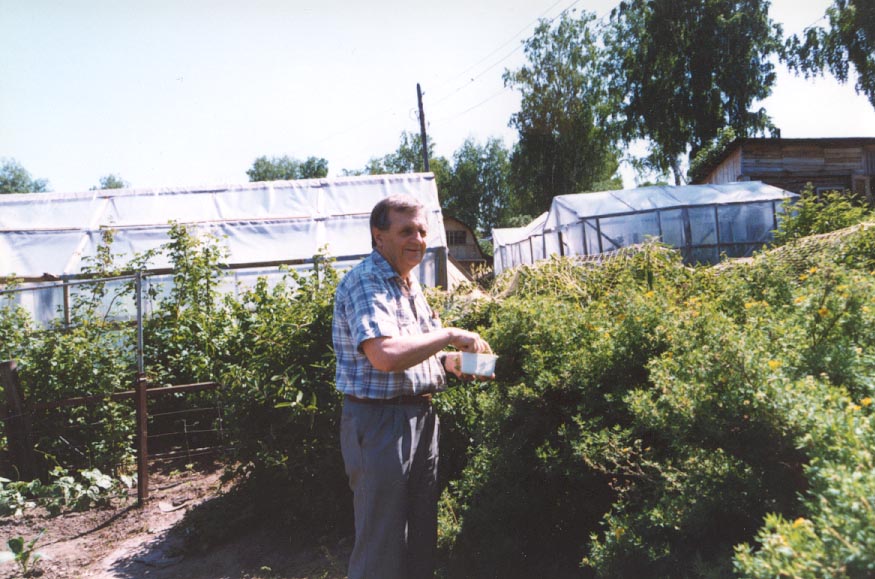}
      }
      \caption{Sergey Godunov collecting fruit at his {\it dacha},  a  typical Russian country seasonal retreat, with gardens 
      and vegetable plots. Novosibirsk,  Siberia, June 1998}
      \label{fig:GodunovDacha}
\end{figure}

For a 70-year-old man, I thought Sergey looked remarkably fit, walking briskly up and down gentle hills. It was a wonderful day, with sunshine and very pleasant temperatures. Upon arrival at his dacha, he introduced me to his family: his wife Tanya, his daughter Ekaterina and husband Alex.  He then showed me around the dacha.  It was a large green plot, abundant with flowers, fruits, and vegetables, laid out with well-organized paths.  He explained that it was his wife who took care of the garden and vegetables. 

Properly processed,  Tanya would preserve much of the fruit and vegetables for the long Siberian winter. He also explained that, traditionally, the dacha in Russia was a real institution—a place for relaxation and perhaps for practising the hobby of cultivating vegetables.  In recent times, he said, it has become a necessity.
I took photographs. We then went into his countryside study, located on a slight elevation of the terrain, with lovely views over the dacha. It was rather hot inside.  

I enjoyed,  as recorded in my diary,  {\it an exquisite and healthy lunch prepared and served with much devotion by Tanya},  largely based on a variety of vegetables.  A very hospitable family,  we all felt at ease and we actually managed to talk a lot. Then it was time to return to Novosibirsk.  Alex gave us a lift in his Lada.
\begin{figure}[h]
      \centerline{
            \includegraphics[scale=0.72, angle=0]{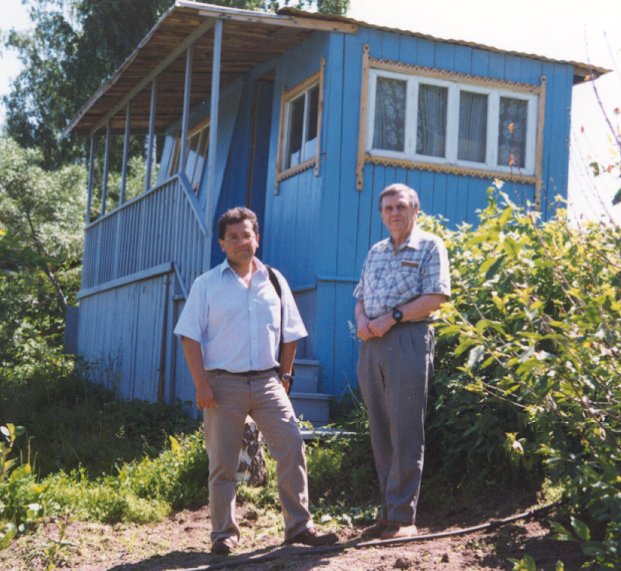}
      }
      \caption{Sergey Godunov (right) and  Eleuterio Toro at Godunov's summer study in his {\it 
      dacha}.   Novosibirsk,  Siberia, June 1998}
      \label{fig:RusanovSchemes}
\end{figure}

Before saying goodbye, we discussed arrangements for two forthcoming events. First,  the celebration of his 70th birthday, to be marked by a special conference to take place in July of the following year in Novosibirsk; he kindly invited me to attend.  On my part, I invited him to visit the UK in September/October of the following year.  As part of his visit, I would also organize a conference in his honour at the University of Oxford. The details were left to be worked out.\\

\begin{center}
 {\bf Godunov's 70th Birthday in Novosibirsk}
\end{center}

In June of 1999, a mathematics conference in Godunov's honour was held at the University of Novosibirsk. At that time, I was spending four months in Japan as a visiting professor at the Tohoku University, Sendai. I encountered some visa-related bureaucratic problems that required me to have a double-entry visa, which I did not have,   to be able to leave Japan and return within a week.  Fortunately, the problem was somehow solved and I flew from Tokyo all the way west to Moscow and then back east to Novosibirsk, on a flight of about four extra hours.

The conference was essentially an all-Russian event.  I identified only two other foreign participants: Rolf Jeltsch (1945–2024) from ETH Zurich and Paul Arminjon (1941–2011) from the University of Montreal.   At the time, Jeltsch was the president of the European Mathematical Society. 

At the conference, Godunov was naturally the centre of attention, attracting scientists, important politicians, the press, and Russian television, which conducted interviews with him. I stood close to one of his TV interviews. From the translation offered by my new Novosibirsk friend Eugene Romenski, Godunov spoke passionately about the loss to his country caused by the exodus of Russian scientists to Western countries. In later years, Godunov would jokingly lament that Toro had “taken away” his collaborators—a reference to long-term visits by Eugene Romenski, first to Manchester and then to Trento, Italy, on several occasions.

There was also time for other social activities. The head of the Sobolev Mathematics Institute of the University of Novosibirsk,  Professor Mikhail Mikhailovich Lavrentiev, invited Rolf Jeltsch, Paul Arminjon,  and myself to meet him in his office.   He welcomed us warmly and offered generous drinks—very nice, I must say, though not exactly the {\it English tea} I had been expecting that warm mid afternoon,  but not less pleasant I must admit.   We all felt quite happy.  

Professor Lavrentiev was the son of mathematician  Michail  Alekseyevich Lavrentiev,  one of the founders of the University of Novosibirsk,  and his son,  named  Michail Michailovich Lavrentiev as well,  was also a mathematician at the Sobolev Mathematics Institute —remarkable, three generations of Lavrentiev mathematicians! The youngest Lavrentiev is currently a Corresponding Member of the Russian Academy of Sciences.

\begin{center}
{\bf Godunov in Manchester,  United Kingdom}
\end{center}

In early September of 1999, Sergey Godunov landed at Manchester Airport to start his one-month visit to the UK. The arrivals screen informed me that his plane had just landed. I rushed to the exit and managed to catch sight of him coming out with his hand luggage—but no suitcase! Controlling my concern, even before greeting and welcoming him to Manchester, I asked him where his luggage was. Smiling, he calmed me down by clarifying that his luggage was the sports bag hanging from his shoulder. 

The idea was to make the most of his visit to the UK. To start with, I organized a local one-day workshop in Manchester itself, with participants mainly from other  universities in the Greater Manchester area. This was followed by short visits to Leeds, Sheffield, London, and Cambridge. His UK visit would be completed with a one-week conference in his honour, to be held at the University of Oxford.
\begin{figure}[h]
      \centerline{
            \includegraphics[scale=0.52, angle=0]{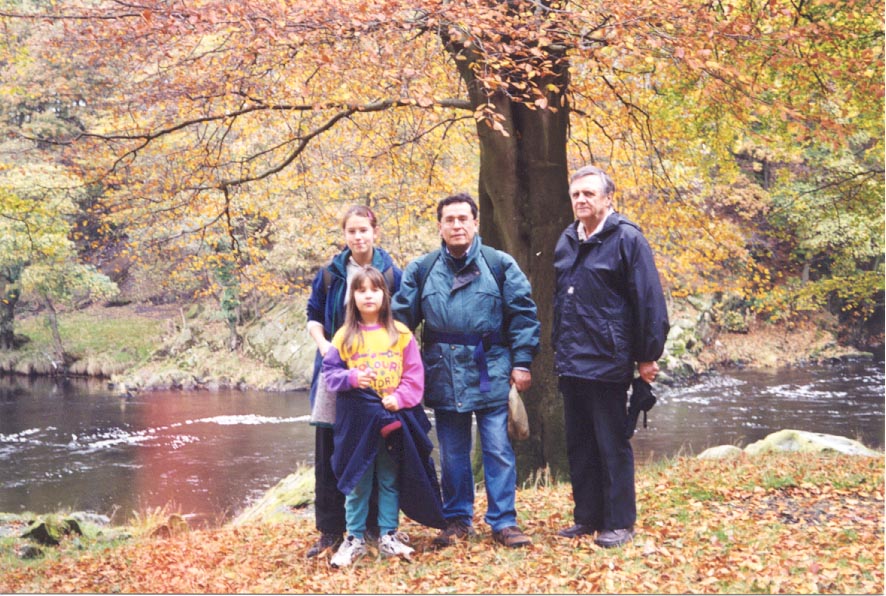}
      }
      \caption{Sergey Godunov (right), Eleuterio Toro and daughters Violeta and Eva.  The Peak District, UK, September 1999}
      \label{fig:GodunovDacha}
\end{figure}
\begin{figure}[h]
      \centerline{
            \includegraphics[scale=0.32, angle=0]{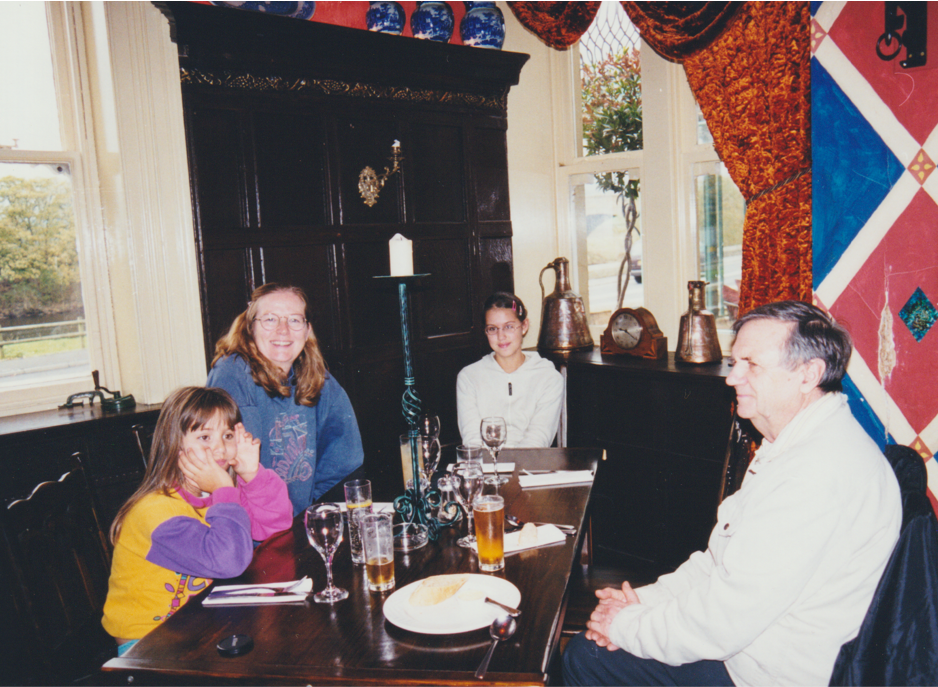}
      }
      \caption{Sergey Godunov, Violeta, Brigitte and Eva.  Pub in the Peak District, UK, September 1999}
      \label{fig:GodunovDacha}
\end{figure}

As part of the social programme I invited Sergey for dinner to my house  one evening.  We also spent a full day out in the Peak District, a characteristic countryside area near Manchester, with my wife Brigitte and my two youngest daughters, Violeta and Eva. The weather was miserable—not rare in Manchester—as it rained all day. The weather of the Siberian summer was surely much better. But Sergey was well prepared:  raincoat, the right kind of shoes for walking in the countryside, and even an umbrella. I always wondered how he managed to have the right attire for the right occasion, all out of his hand luggage.

Enjoying a few pints of bitter in a typical English pub was an experience that Professor Godunov could not miss, I thought.  Some colleagues and friends joined in  for a drink, one evening. The brick fireplace next to our table, even if not lit at that time of the year, somehow gave the place a cosy atmosphere.  Sergey told us about his Moscow years as a mathematics student,  and about his teachers,  such as
I.G.  Petrovsky,  B.N. Delaunay, I.M. Gelfand,and others.   His PhD examiner was S. L.  Sobolev.

As an undergraduate mathematics student,  at the age of 19,  he published his first paper on a number theory subject,  jointly with his teacher B.N.  Delaunay.   As a researcher he interacted closely with renown scientists, such as D. A.  Frank-Kamenetskii and Y. B.  Zel'dovich, amongst others.

As time went on that evening,  we all relaxed a bit and Sergey’s English kept flowing better and better.  He also spoke about the early years of his career at the Novosibirsk University,  which he started in 1969.  He told us that he and his friend  Boris Delaunay,  were great enthusiasts of the Russian countryside and of long walks in rugged,  mountainous terrain.  Every week they would take a couple of days off, sleeping rough when night arrived.  He told us of an accident he had suffered while hiking in the Caucasus Mountains: he fell from a rocky elevation, suffering a serious head injury,  the signs of which were still evident—a deep,  unnatural hollow in his head. 

He allowed us, one after the other, to feel his head and find the indentation.  Some local customers kept giving us funny looks, perhaps wondering what was so interesting about this man’s head.   These Mancunians did not know, of course, that the man's head was no ordinary head—it sheltered a beautiful mind.  Godunov survived that accident. There was a brief silent moment. Surely, we all felt relieved,  perhaps reflecting on  how fortunate mankind was.

It may sound strange, but it was not easy to find technical issues in common to discuss with Sergey. He was not particularly interested in his old work on Godunov’s method.  Rather,  he was engaged in new subjects, such as grid generation and computational linear algebra. He was particularly interested in moving meshes to capture (track) discontinuities with first-order methods. I myself was busy developing what we came to call ADER methods  started in the early 1990s but first communicated in the early 2000s \cite{Toro:2001c, Toro:2002a, Titarev:2002a}.  

I attempted to explain that ADER methods were a perfect higher-order analogue  of Godunov’s first-order method.  I explained that such methods retained the two distinguishing features of his method: the update formula,  called {\it Numerical method} in Fig.  \ref{fig:GodunovMethod},  was still a one-step formula,  and the classical piecewise-constant-data {\it Riemann problem} of his method was replaced by a generalized Riemann problem (GRP$m$), with initial data consisting of polynomials of degree $m$, or some other suitable functions.  I also dared to point out that the GRP$m$  included source terms,  if present in the original equations.

Encouraged by his attentiveness and politeness when listening I went on to say that  available results for linear problems were very promising,  but non-linear extensions were proving challenging.  At that stage I sensed,  however,  that Sergey was not particularly impressed with my description of a very-high-order extension of Godunov's method.   We left it at that on that occasion and went for a cup tea and biscuits.

\begin{center}
{\bf Godunov's 70th Birthday in Oxford,  UK}
\end{center}

At St.  Anne’s College, Oxford University, UK, 140 scientists from 20 countries gathered to commemorate Godunov’s 70th birthday. It was a grand occasion; one enthusiastic participant even commented, “Everyone is here!” Regrettably, however, there were some important absentees. 

The conference took place from 18th to 22nd October 1999. Plenary speakers included S. K. Godunov, R. Abgrall, H. Aiso, P. Colella, D. Drikakis, P. García-Navarro, B. Gustafsson, B. Koren, D. Kroener, A. Marquina, S. Osher, P. A. Raviart, P. L. Roe, E. Romenski, T. Saito, R. Saurel, P. K. Sweby, E. Tadmor, E. F. Toro, I. Toumi, L. N. Trefethen, M. E. Vázquez, and B. Wendroff.
The conference programme was supplemented by many contributed high-level talks and was preceded by an intensive three-day short course on topics related to the Godunov conference, and was well attended,  predominantly by young researchers,  but also by senior academics.

\begin{figure}[h]
      \centerline{
               \includegraphics[scale=0.39, angle=0]{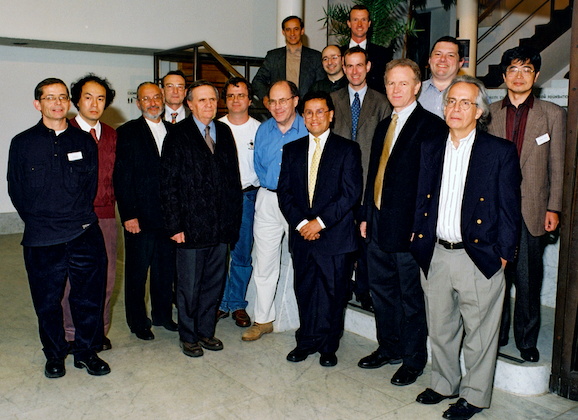}
      }
      \caption{Sergey Godunov (fifth from left) and speakers at the {\it International conference in honour of S.  K.  Godunov on occasion of his 70th birthday}.  Oxford University,  UK, October 1999}
      \label{fig:RusanovSchemes}
\end{figure}
The main themes of the conference were related to Godunov’s 1959 seminal paper, which introduced what became known as Godunov’s numerical method for solving hyperbolic equations \cite{Godunov:1959a}.  During the week-long event, it became increasingly clear that Sergey K. Godunov had also made very significant contributions in other areas of mathematics.  One of the plenary speakers remarked, “I know no other living scientist who has made so many contributions to so diverse areas of research.” 

Godunov's pioneering work on the theory of hyperbolic partial differential equations,  thought less known than Godunov's numerical method,  was also highlighted through discussions and two relevant contributions to the conference,  one by Godunov himself \cite{Godunov:2001a} and another by his close collaborator, and first PhD student,  Eugene Romenski  \cite{Romenski:2001a}. 

Back in 1961 Godunov proposed classes of systems of partial differential equations for which he  developed a corresponding ground-braking  mathematical theory \cite{Godunov:1961a}.  These theoretical advances emanated from his previous studies on the relationship between thermodynamics and the well posedness of the equations of mathematical physics.  His work in this direction actually started from his seminal paper of 1959 \cite{Godunov:1959a}  that gave rise to the famous Godunov’s method for computing solutions to non-linear hyperbolic equations,  as already mentioned.   

Remarkably,  his 1961 theory \cite{Godunov:1961a} anticipated by ten years subsequent communications on the same subject  \cite{Friedrichs:1971a}; see  also \cite{Busto:2022a, Warnecke:2025a}.  Later developments were carried out in collaboration with Novosibirsk researchers Eugene Romenski,  Sergey Gavrilyuk and Godunov's last PhD student Ilya Peshkov; see for example \cite{Godunov:1972a, Peshkov:2016a, Chiocchetti:2021a}.  Of course,  further developments of Godunov's seminal works have been carried out by other groups,  see  \cite{Ruggeri:1981a},  for example.   A very recent review on some of these developments is found in \cite{Warnecke:2025a}.

The Oxford conference gave rise to a 1077-page compendium on {\it Godunov Methods: Theory and Applications}, which included 97 refereed papers arising from the conference presentations \cite{Toro:2001b}. The event was supported, among others, by the London Mathematical Society.

\begin{center}
{\bf Godunov at the Newton Institute,  \\
University of Cambridge,  UK}
\end{center}

In 2003,  C. M. Dafermos (Brown University, USA), P. G. LeFloch (Centre de Physique Théorique,  de l'École Polytechnique, France), and E. F. Toro (University of Trento, Italy),  organized a six-month research programme at the Newton Institute for Mathematical Sciences, University of Cambridge, UK, under the title {\it Nonlinear Hyperbolic Waves in Phase Dynamics and Astrophysics}. This was a unique occasion that brought  together world leaders from various scientific fields related to the programme. More information can be found at https://www.newton.ac.uk/event/npa/
\begin{figure}[h]
      \centerline{
              \includegraphics[scale=0.37, angle=0]{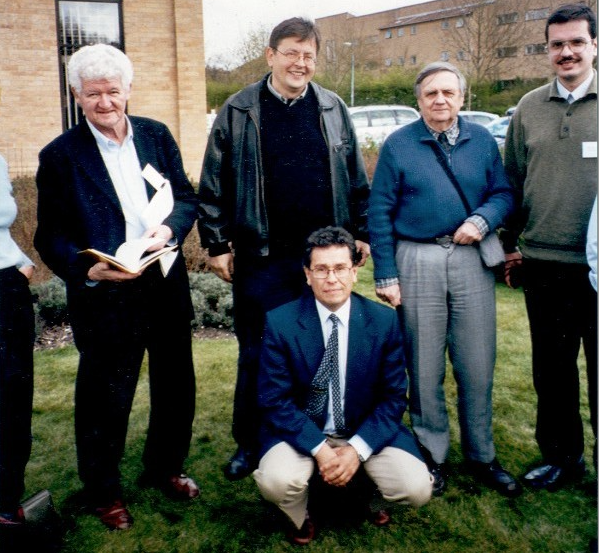}
      }
      \caption{Sergey Godunov (second from right),  Peter Lax (first from left), Vladimir Titarev (first from right) and Eleuterio Toro (squatting).  Research programme  
      {\it  Nonlinear hyperbolic waves in phase dynamics and astrophysics,  organized by C.  M. Dafermos,  P.  G.  LeFloch and E. F.   
      Toro. }
      Newton Institute for Mathematical Sciences, University of Cambridge, UK, 2003}
      \label{fig:Godunov-Lax}
\end{figure}

Sergey Godunov was one of our key participants, as were Peter Lax (1926-2025), Alberto Bressan, Barbara Lee Keyfitz,  Andrea Corli,  Randall LeVeque, Chi-Wang Shu,  
Peter Sweby, Ash Kapila, John Bell, Blake Temple, Shi Jin, Eugene Romenski, Marshall Slemrod, Gui-Qiang Chen,  Tai-Ping Liu, and many others. The programme also included a special lecture by the late Stephen Hawkins (1942-2018).

There was much anticipation in the air at the prospect of seeing two giants together—Sergey Godunov and Peter Lax—especially in light of past controversies.  It was a historic moment.  See photo \ref{fig:Godunov-Lax}.
\vspace{-4mm}
\begin{center}
{\bf Godunov's 90th Birthday,  Novosibirsk}
\end{center}

In August 2019,  more than 300 scientists from 24 countries gathered at a major international conference in Novosibirsk to honour S. K. Godunov,  on occasion of his 90th birthday.   It was an extraordinary event,  with Godunov himself giving the one-hour inaugural lecture.  He spoke about {\it Memories of finite difference schemes}.  I had the honour of giving one of the plenary lectures,  titled 
{\it Godunov's methods: the ADER path to high-order of accuracy} \cite{Toro:2020a}.  Of all the presentations given,  52 were selected for publication in the book proceedings of the conference,  titled {\it Continuum Mechanics, Applied Mathematics, Scientific Computing: Godunov's Legacy},  edited by G. V. Demidenko, E. Romenski, E. Toro and M. Dumbser (Springer, 2020) \cite{Demidenko:2020a}.  

This was the last time I saw Sergey Godunov.  He passed away on 15th July 2023.  Curiously,  his 94th birthday was due two days later, on 17th July.
\begin{figure}[h]
      \centerline{
           \includegraphics[scale=0.154, angle=0]{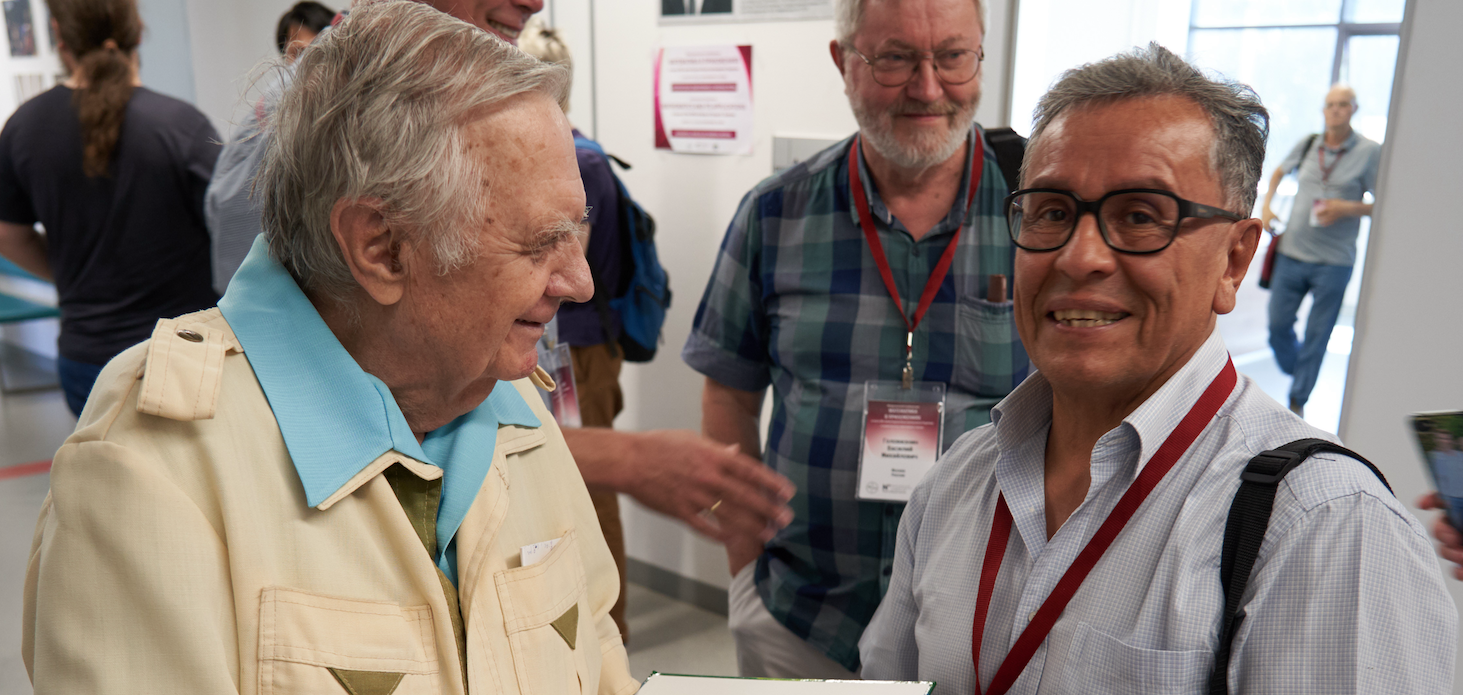}
      }
      \caption{S. K. Godunov and E. F. Toro at the international conference  to honour  Godunov's  90th birthday.   Novosibirsk (Siberia, Russia),  
      4th-10th August 2019}
      \label{fig:Godunov90}
\end{figure}
\begin{figure}[h]
      \centerline{
           \includegraphics[scale=0.2, angle=0]{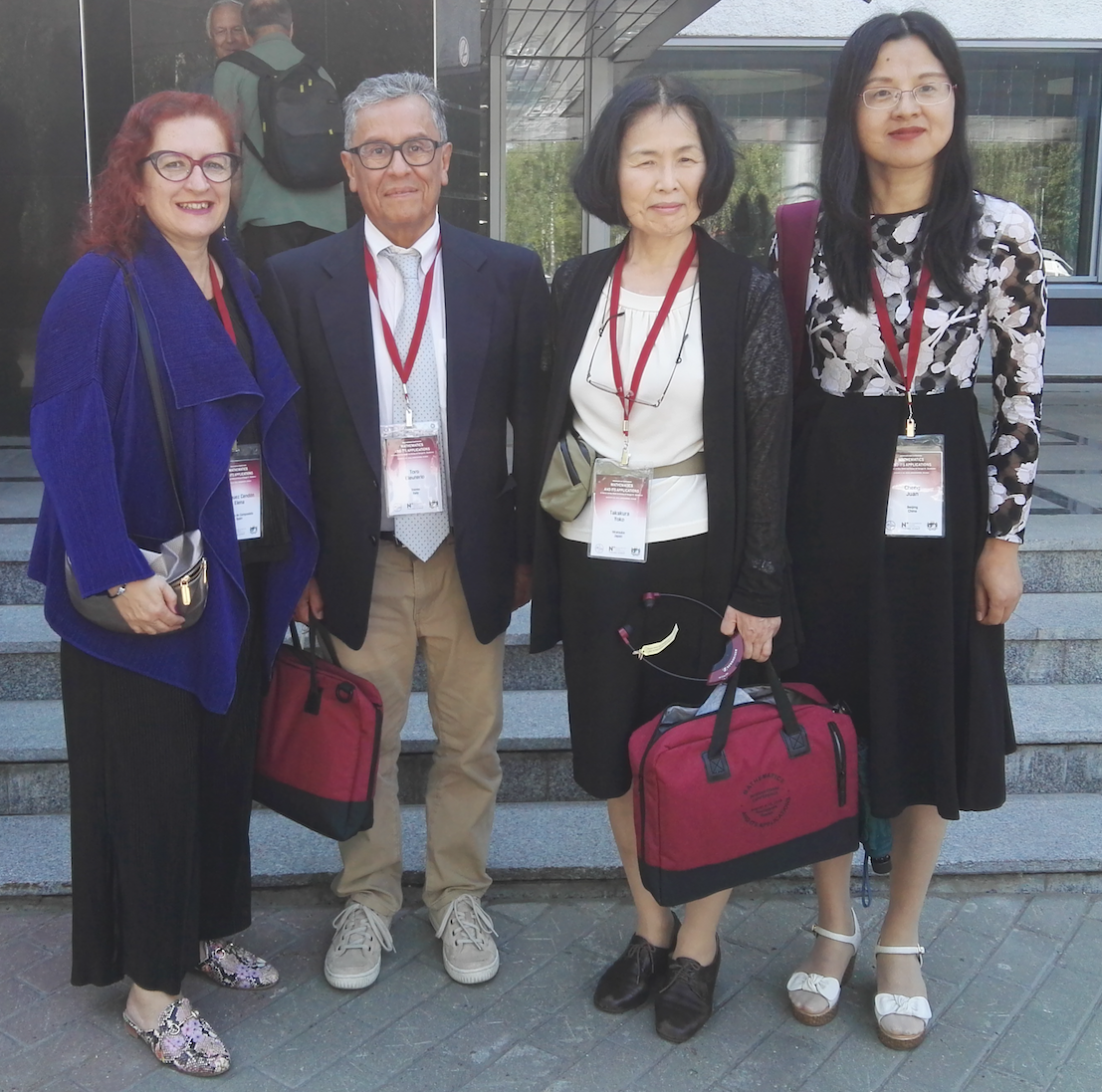}
      }
      \caption{From right: Juan Cheng (China); Yoko Takakura (Japan); Eleuterio Toro (Italy) and Elena V\'azquez (Spain).  International 
      conference  to honour  Godunov's  90th birthday.   Novosibirsk (Siberia, Russia),  4th-10th August 2019}
      \label{fig:ChengEtal}
\end{figure}
\vspace{-5mm}
\begin{center}
{\bf Novosibirsk-Manchester-Trento Collaboration}
\end{center}

During my first meetings with Godunov and his collaborators  in Novosibirsk (1998,1999), we discussed potential collaboration and visitor exchanges.   Through funding from a European project  I made arrangements for a visiting position in Manchester for Eugene Romenski,  a professor at the 
Sobolev Institute of Mathematics in Novosibirsk.   As already noted,  Romenski was Godunov's first PhD student in Novosibirsk;  he defended his thesis in 1975.   Eugene  took a temporary position as  Senior Research Fellow,  at the Manchester Metropolitan University,  UK,  from April 2001  to October 2002.   After I moved to the University of Trento,  Italy,  in 2002,  I continued the collaboration with Romenski.  He held  a visiting professorship in my Department from September 2003  to October 2005. This project was funded by the Italian Ministry of Education, University  and Research.  We worked on mathematical models for two phase compressible flows following the Godunov framework of symmetric hyperbolic systems \cite{Romenski:2004a, Romenski:2007a}.

The arrival of Michael Dumbser from Stuttgart to my group in Trento in 2006,  first as a post-doctoral fellow,  then as a research assistant  and later as a full professor,  gave increased momentum to the collaboration between the Novosibirsk school and our Trento group.   Importantly,  this phase departed  from the early works of Godunov and Romenski on sophisticated, unified,  mathematical models for fluid and solid mechanics.   A further impulse to these developments were given by the arrival  of Ilya Peshkov to the Trento group,  on a permanent position as associate professor.  I note that Ilya was Godunov's last PhD student at the Sobolev Mathematics Institute in Novosibirsk (2007-2009). 

Significantly,  the deployment of advanced numerical methodology developed at Trento,   to accurately solve the equations of the sophisticated models of Godunov,  Romenski and Peshkov has given rise to a very substantial area of research,  with important scientific advances,  all largely led by my colleague Michael Dumbser \cite{Busto:2022a,Dumbser:2016b, Dumbser:2016c}.  

\begin{figure}[h]
      \centerline{
           \includegraphics[scale=0.2, angle=0]{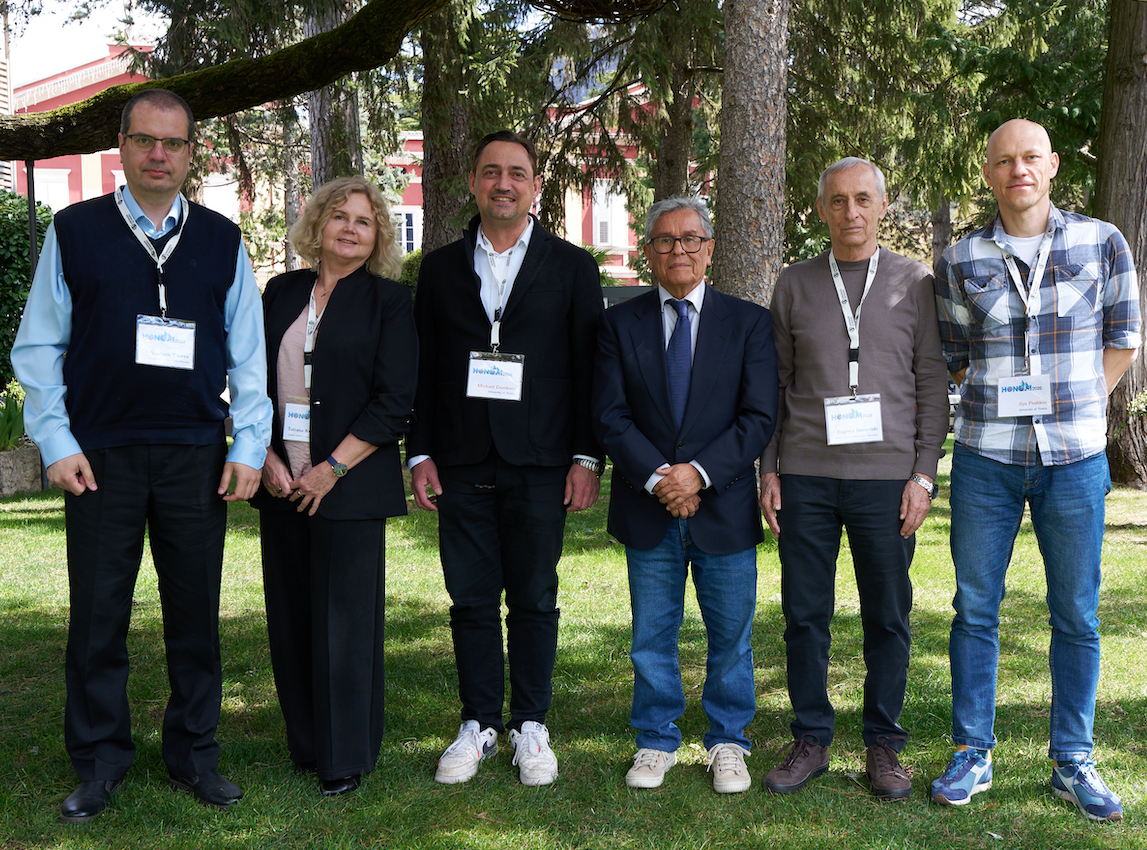}
      }
      \caption{From right to left: Ilya Peshkov, Eugene Romenski, Eleuterio Toro, Michael Dumbser, Tatiana Kozubskaya and Vladimir Titarev.  International Conference on High Order Nonlinear Numerical Methods for Evolutionary PDE- HONOM 2026, 30th March to 4th April 2026. Trento, Italy.}
      \label{fig:GodunovSpeaking}
\end{figure}

\begin{center}
{\bf Concluding Reflections}
\end{center}

The Russian mathematician Sergey K. Godunov made profound contributions to science, spanning physics, partial differential equations,  numerical analysis and scientific computing. His work established a research paradigm that, over the past six decades, has inspired distinctive lines of inquiry pursued by researchers around the world. As a university academic, he trained generations of scientists who lived up to the high standards he set.

As is the case with many groundbreaking scientific discoveries, some of Godunov’s contributions have been accompanied by misunderstandings, prejudice, and even controversy. Godunov’s method continues to encounter resistance within parts of the numerical analysis community. In particular, the Godunov numerical flux is often narrowly interpreted as relying on the exact solution of the Riemann problem for its evaluation. This interpretation has led to the widespread belief that Godunov’s method is excessively complicated, computationally expensive, and impractical for demanding applications in which an exact Riemann solver may be unavailable or prohibitively costly.

This narrow view on Godunov's method, however, overlooks the fact that in 1962—only three years after the original publication of his method—Godunov was also the first to derive and communicate a linearized approximate Riemann solver.   Moreover, his fundamental contributions to the mathematical theory of hyperbolic partial differential equations remained unrecognized for many years, as documented in the literature \cite{Godunov:1961a, Friedrichs:1971a}; see also \cite{Busto:2022a, Warnecke:2025a,Ruggeri:2025a}.

Language barriers are partly to blame for this, as they are detrimental to all forms of human expression. This was already well recognized in the 1970s by mathematician Sir Christopher Zeeman at the Warwick Mathematics Institute, particularly in the context of mathematics education. Geopolitics has likewise been a persistent obstacle to openness in the sharing of human creativity across the arts, sports, science, and culture more broadly. Godunov’s entire life unfolded during an especially turbulent period of modern history, marked by profound divisions within humanity.

Nonetheless, the power of Godunov’s creativity has influenced generations of scientists around the world, despite these barriers. His ideas have shaped the lives of many individuals in various ways, ranging from university courses to the choice of a final-year project, a master’s thesis, or a PhD dissertation—and beyond.  

Some have pursued careers in industry, where mathematical modelling and computation are used to solve practical engineering problems, and where Godunov’s PDE theory and Godunov-type computational methods may play a crucial role.

Others have pursued academic careers, ranging from postdoctoral research fellowships to lectureships and professorships. Each, in turn, has acted as a multiplying force, expanding an ever-growing network of scientists working on the theory of partial differential equations, related numerical analysis, and their applications—all inspired by Godunov’s original ideas.

At a fundamental level, Sergey K. Godunov gave many of us a dignified occupation for life.

\begin{center}
{\bf Acknowledgements}
\end{center}

I am indebted to Violeta Da Rold, Eugene Romenski, Annunziato Siviglia and Tommaso Ruggeri for their help in completing this document.

\bibliographystyle{plain}
\bibliography{refs-latest}

\end{document}